\documentstyle[11pt,amssymb]{article}
\begin{document}

{\large

\noindent{\bf Big $q$-Laguerre and $q$-Meixner  polynomials and
representations}

\noindent{\bf of the algebra $U_q({\rm su}_{1,1})$}  }
\bigskip

{\bf M. N. Atakishiyev, N. M. Atakishiyev, and A. U. Klimyk}

Instituto de Matem\'aticas, UNAM, CP 62210 Cuernavaca, Morelos,
M\'exico

\bigskip

E-mail: natig@matcuer.unam.mx and anatoliy@matcuer.unam.mx
\bigskip
\bigskip

\begin{abstract}
Diagonalization of a certain operator in irreducible
representations of the positive discrete series of the quantum
algebra $U_q({\rm su}_{1,1})$ is studied. Spectrum and
eigenfunctions of this operator are found in an explicit form.
These eigenfunctions, when normalized, constitute an orthonormal
basis in the representation space. The initial $U_q({\rm
su}_{1,1})$-basis and the basis of eigenfunctions are interrelated
by a matrix with entries, expressed in terms of big $q$-Laguerre
polynomials. The unitarity of this connection matrix leads to an
orthogonal system of functions, which are dual with respect to big
$q$-Laguerre polynomials. This system of functions consists of two
separate sets of functions, which can be expressed in terms of
$q$-Meixner polynomials $M_n(x; b,c;q)$ either with positive or
negative values of the parameter $b$. The orthogonality property
of these two sets of functions follows directly from the unitarity
of the connection matrix. As a consequence, one obtains an
orthogonality relation for the $q$-Meixner polynomials
$M_n(x;b,c;q)$ with $b<0$. A biorthogonal system of functions
(with respect to the scalar product in the representation space)
is also derived.

\end{abstract}


\bigskip
\bigskip
\bigskip

\noindent{\bf 1. Introduction}
\bigskip

The significance of representations of Lie groups and Lie algebras
for studying orthogonal polynomials and special functions is well
known. The appearance of quantum groups and quantized universal
enveloping algebras (quantum algebras) and development of their
representations led to their applications for elucidating
properties of $q$-orthogonal polynomials and $q$-special functions
(see, for example, [1--3]). Since the theory of quantum groups and
their representations is much more complicated than the Lie
theory, the corresponding applications in the theory of orthogonal
polynomials and special functions are more difficult. At the first
stage of such applications, the compact quantum groups and their
finite dimensional representations have been used.

It is known that the noncompact Lie group $SU(1,1)\sim SL(2,{\Bbb
R})$ and its representations are very productive for the theory of
orthogonal polynomials and special functions. Unfortunately, there
are difficulties with a satisfactory definition of the noncompact
quantum group (or algebras of functions on the quantum group)
$SU_q(1,1)$, which would give us a possibility to use such quantum
group extensively for deeper understanding the theory of
orthogonal polynomials and special functions. For this reason,
representations of the quantum algebra $U_q({\rm su}_{1,1})$ have
been commonly used for such purposes (see, for example, [4--7]).

In this paper we use representations of the positive discrete
series of the quantum algebra $U_q({\rm su}_{1,1})$ for exploring
properties of big $q$-Laguerre polynomials and $q$-Meixner
polynomials. In fact, we deal with certain operators in these
representations and do not touch the Hopf structure of the algebra
$U_q({\rm su}_{1,1})$. Our study of these polynomials is related
to representation operators, which can be represented by a Jacobi
matrix. We deal with those properties of big $q$-Laguerre
polynomials and $q$-Meixner polynomials, which have not been
usually investigated by means of group theoretical methods.
Namely, we consider a representation operator $A$, represented by
some particular Jacobi matrix. We diagonalize this selfadjoint
bounded operator with the aid of big $q$-Laguerre polynomials. The
orthogonality relation for these polynomials allows us to find a
spectrum of the operator $A$. It is simple and discrete. We find
an explicit form of all eigenfunctions of this operator. Since the
spectrum is simple, these eigenfunctions constitute an orthogonal
basis in the representation space. Then we normalize this basis.
As a result, one has two orthonormal bases in the representation
space: the canonical (or the initial) basis and the basis of
eigenfunctions of the operator $A$. They are interrelated by a
unitary matrix $U$ with entries $u_{mn}$, which are explicitly
expressed in terms of big $q$-Laguerre polynomials. Since
$U^*\,U=U\,U^*=I$, there are two requirements for providing the
unitarity of this matrix $U$, namely:
$$
\sum _{n=0}^\infty u_{mn}u_{m'n}=\delta_{mm'},\ \ \
\sum _{m=0}^\infty u_{mn}u_{mn'}=\delta_{nn'}.   \eqno (1.1)
$$
The first relation expresses the orthogonality relation for big
$q$-Laguerre polynomials. In order to interpret the second
relation, we consider big $q$-Laguerre polynomials
$P_n(q^{-m};a,b;q)$ as functions of $n$. In this way one obtains a
set of orthogonal functions, which are expressed in terms of two
sets of $q$-Meixner polynomials (these two sets of $q$-Meixner
polynomials can be considered as a dual set of polynomials with
respect to big $q$-Laguerre polynomials: such duality property is
well known in the case of discrete polynomials, orthogonal on a
finite set of points). Consequently, the second relation in (1.1)
leads to the orthogonality relations for those $q$-Meixner
polynomials, which enter to our two sets (with respect to certain
measures). For one of these sets we obtain the well-known
orthogonality relation for $q$-Meixner polynomials. The second set
leads to an orthogonality relation for $q$-Meixner polynomials
$M_n(x;b,c;q)$ with $b<0$. As far as we know, the possibility of
extending the common orthogonality relation for $q$-Meixner
polynomials $M_n(x;b,c;q)$ to a wider range of the parameter $b$
has not been discussed in the literature.

We deduce from orthogonality relations for  functions, dual to the
big $q$-Laguerre polynomials, that $q$-Meixner polynomials are
associated with the indeterminate moment problem and the commonly
known orthogonality measure for these polynomials is not an
extremal one. (Note that if this measure would be extremal, then
the set of the $q$-Meixner polynomials would form a basis in the
corresponding Hilbert space.)

We consider also the case of two representation operators $A_1$
and $A_2$, such that $A^*_1=A_2$. Diagonalization of these
operators leads to two sets of functions, which are biorthogonal
with respect to the scalar product in the representation space.

Throughout the sequel we always assume that $q$ is a fixed
positive number such that $q<1$. We use the theory of special
functions and notations of the standard $q$-analysis (see,
for example, [8] and [9]). We use $q$-numbers $[a]_q$ defined as
$$
[a]_q = \frac{q^{a/2}-q^{-a/2}}{q^{1/2}-q^{-1/2}}, \eqno (1.2)
$$
where $a$ is any complex number. We shall also use the well known
notations
$$
(a;q)_n:= (1-a)(1-aq)\cdots (1-aq^{n-1}),\ \ \ (a_1,\cdots
,a_k;q)_n=(a_1;q)_n\cdots (a_k;q)_n.
$$
\medskip

\noindent{\bf 2. Discrete series representations of
$U_q({\rm su}_{1,1})$}
\bigskip

 The quantum algebra $U_q({\rm
su}_{1,1})$ is defined as the associative algebra, generated by
the elements $J_+$, $J_-$, and $J_0$, satisfying the
commutation relations
$$ [J_0,J_{\pm}]=\pm J_{\pm},\ \ \ \ [J_-,J_+] =
{q^{J_0}-q^{-J_0}\over q^{1/2}-q^{-1/2}} \equiv [2J_0]_q ,
$$
and the conjugation relations
$$
J_0^* = J_0,\ \ \ \ J_+^* = J_-. \eqno (2.1)
$$
(Observe that here we have replaced $J_-$ by $-J_-$ in the usual
definition of the algebra $U_q({\rm sl}_{2})$.)

We are interested in the discrete series representations of
$U_q({\rm su}_{1,1})$ with lowest weights. These irreducible
representations will be denoted by $T^+_l$, where $l$ is a lowest
weight, which can take any positive number (see, for example,
[10]).

The representation $T^+_{l}$ can be
realized on the space ${\cal L}_{l}$ of all polynomials in $x$. We
choose a basis for this space, consisting of the monomials
$$
f^l_n\equiv f^{l}_n(x) := c^{l}_n\, x^n, \ \ \  n = 0,1,2,\cdots ,
\eqno(2.2)
$$
where
$$ c^{l}_0 = 1, \qquad c^{l}_n = \prod_{k=1}^n\,{[2l+k-1]_q^{1/2}
\over [k]_q^{1/2}} = q^{(1-2l)n/4}{(q^{2l};q)_n^{1/2}\over
(q;q)_n^{1/2}}\,, \ \ n = 1,2,3, \cdots .     \eqno(2.3)
$$
The representation $T^+_{l}$ is then realized by the operators
$$
J_0 =x{d\over dx} + l,\qquad J_{\pm} = x^{\pm 1}[ J_0(x) \pm l]_q
\,.
$$
As a result of this realization, we have
$$
J_+\, f^{l}_n =\sqrt{[2l+n]_q \,[n+1]_q }\, f^{l}_{n+1}
=\frac{q^{-(n+l-1/2)/2}}{1-q} \sqrt{(1-q^{n+1})(1-q^{2l+n})}
f^{l}_{n+1} , \eqno (2.4)
$$ $$
J_-\, f^{l}_n =\sqrt{[2l+n-1]_q \,[n]_q }\, f^{l}_{n-1}
=\frac{q^{-(n+l-3/2)/2}}{1-q} \sqrt{(1-q^{n})(1-q^{2l+n-1})}
f^{l}_{n-1},  \eqno (2.5)
$$  $$
J_0\, f^{l}_n = (l + n)\,f^{l}_n . \eqno (2.6)
$$

We know that the discrete series representations $T_l$ can be
realized on a Hilbert space, on which the conjugation relations
(2.1) are satisfied. In order to obtain such a Hilbert space, we
assume that the monomials $f^{l}_n(x)$, $n=0,1,2,\cdots$,
constitute an orthonormal basis for this Hilbert space. This
introduces a scalar product $\langle \cdot ,\cdot \rangle$ into
the space ${\cal L}_l$. Then we close this space with respect to
this scalar product and obtain the Hilbert space, which will be
denoted by ${\cal H}_l$.
\bigskip

\noindent{\bf 3. Representation operators related to
big $q$-Laguerre polynomials}
\bigskip

In this section we are interested in the operator
$$
A :=\alpha q^{J_0/4}\,\left(\sqrt{1-bq^{J_0-l}}J_+ \,q^{(J_0-l)/2} +
q^{(J_0-l)/2} J_- \sqrt{1-bq^{J_0-l}}\,\right)\, q^{J_0/4}
- \beta_1 q^{2J_0} + \beta_2 q^{J_0-l}  \eqno (3.1)
$$
of the representation $T^+_l$, where $b<0$ and
$$
\alpha =(-b)^{1/2}q^{l}(1-q),\ \ \
\beta_1=b(1+q),\ \ \
\beta_2=bq+q^{2l}(b+1).
$$
Since the bounded operator $q^{J_0}$ is diagonal in the basis
$f^l_n$, $n=0,1,2,\cdots$, without zero diagonal elements, the
operator $A$ is well defined.

We have the following expression for the symmetric operator $A$ in
the canonical basis $f^l_n$, $n=0,1,2,\cdots$:
$$
A\, f^l_n=(-ab)^{1/2}q^{(n+2)/2}\biggl[
\sqrt{(1-q^{n+1})(1-aq^{n+1})(1-bq^{n+1})} f^{l}_{n+1}
$$  $$
+\, q^{-1/2}\sqrt{(1-q^{n})(1-aq^{n})(1-bq^n)} f^{l}_{n-1}\biggr]
$$  $$
-[ abq^{2n+1}(1+q)-q^{n+1}(a+ab+b)] f^l_n, \ \ \ \ a=q^{2l-1}.
$$
Since $q<1$ the operator $A$ is bounded.
Therefore, one can close this operator and we assume in what
follows that $A$ is a closed
(and consequently defined on the whole space ${\cal H}_l$)
operator. Since $A$ is symmetric, its closure is a selfadjoint
operator.

We wish to find eigenfunctions $\xi_\lambda (x)$ of the operator
$A$, $A \xi_\lambda (x)=\lambda \xi_\lambda (x)$. We set
$$
\xi_\lambda (x)=\sum _{n=0}^\infty a_n(\lambda)f^l_n (x).
$$
Acting by the operator $A$ upon both sides of this relation,
one derives that
$$
\sum _{n=0}^\infty a_n (\lambda)\Bigl\{ q^{(n+2)/2}(-ab)^{1/2}
\sqrt{(1-q^{n+1}) (1-aq^{n+1})(1-bq^{n+1})}f^l_{n+1}
$$  $$
+\, q^{(n+1)/2}(-ab)^{1/2}
\sqrt{(1-q^{n}) (1-aq^{n})(1-bq^n)}f^l_{n-1}
+d_nf^l_n\Bigr\} =\lambda \sum _{n=0}^\infty
a_n(\lambda)f^l_n,
$$
where
$$ d_n=-abq^{2n+1}(1+q)+ q^{n+1}(a+ab+b).
$$
Note that $0<a<q^{-1}$ since $l$ can be any positive number.
Comparing coefficients of a fixed $f^l_n$, one obtains a
three-term recurrence relation for the coefficients $a_n
(\lambda)$:
$$
q^{(n+2)/2}(-ab)^{1/2}
\sqrt{(1-q^{n+1}) (1-aq^{n+1})(1-bq^{n+1})}
a_{n+1}(\lambda)
$$ $$
+\, q^{(n+1)/2}(-ab)^{1/2}
\sqrt{(1-q^{n})(1-aq^{n})(1-bq^{n})}
a_{n-1}(\lambda)+d_na_n (\lambda)=\lambda a_n (\lambda).
$$
We make here the substitution
$$
a_n (\lambda)=(-ab)^{-n/2}q^{-n(n+3)/4} \left( \frac{(aq,bq;q)_n}{(q;q)_n }
\right) ^{1/2} a'_n (\lambda)
$$
and derive the relation
$$
(1-aq^{n+1})(1-bq^{n+1}) a'_{n+1}(\lambda)- abq^{n+1}(1-q^{n})
a'_{n-1}(\lambda) +d_na'_n(\lambda)=\lambda a'_n (\lambda),
$$
where $0<a<q^{-1}$ and $b<0$. It coincides with the recurrence
relation for the big $q$-Laguerre polynomials
$$
P_n(\lambda ;a,b;q):= {}_3\phi_2 (q^{-n},0,\lambda;\;  aq,bq;\;q,q)
$$  $$
=(q^{-n}/b;q)^{-1}_n\; {}_2\phi_1 (q^{-n},\, aq/\lambda ;\; aq;\
q,\lambda /b)   \eqno (3.2)
$$
(see formula (3.11.3) in [11]), that is, $
a'_n(\lambda)=P_n(\lambda;a,b;q)$, $a=q^{2l-1}$. Therefore,
$$
a_n(\lambda)=(-ab)^{-n/2}q^{-n(n+3)/4} \left( \frac{(aq,bq;q)_n}
{(q;q)_n} \right) ^{1/2} P_n( \lambda ;a,b;q) .  \eqno (3.3)
$$

Thus, the eigenfunctions of the operator $A$ have the form
$$
\xi_\lambda(x)= \sum _{n=0}^\infty (-ab)^{-n/2}q^{-n(n+3)/4}\left(
\frac{(aq,bq;q)_n}{(q;q)_n} \right) ^{1/2} P_n( \lambda ;a,b;q)
f^l_n (x)  \eqno (3.4)
$$  $$
=\sum _{n=0}^\infty (-b)^{-n/2}a^{-3n/4}q^{-n(n+3)/4}
\frac{(aq;q)_n} {(q;q)_n} (bq;q)^{1/2}_n P_n( \lambda ;a,b;q)x^n.
$$

To find a spectrum of the operator $A$, we take into account the
following. The selfadjoint operator $A$ is represented by a Jacobi
matrix in the basis $f^l_n(x)$, $n=0,1,2,\cdots$. According to the
theory of operators of such type (see, for example, [12], Chapter
VII), eigenfunctions $\xi_\lambda$ of this operator are expanded
into series in the monomials $f^l_n(x)$, $n=0,1,2,\cdots$, with
coefficients, which are polynomials in $\lambda$. These
polynomials are orthogonal with respect to some measure $d\mu
(\lambda)$ (moreover, for selfadjoint operators this measure is
unique). The set (a subset of ${\Bbb R}$), on which the
polynomials are orthogonal, coincides with the spectrum of the
operator under consideration and the spectrum is simple. Let us
apply these assertions to the operator $A$.

The orthogonality relation for the big $q$-Laguerre polynomials
$P_n( \lambda ;a,b;q)$ for $0<a<q^{-1}$ and $b<0$ is known to be
of the form
$$
\int_{bq}^{aq} \frac{(x/a,x/b;q)_\infty}{(x;q)_\infty} P_m(
x ;a,b;q)P_{m'}( x ;a,b;q)d_qx
$$  $$
\equiv  \sum _{n=0}^\infty \frac{(q^{n+1};q)_\infty
(aq^{n+1}/b;q)_\infty} {(aq^{n+1};q)_\infty}q^n P_m( aq^{n+1}
;a,b;q)P_{m'}( aq^{n+1} ;a,b;q)
$$  $$
-\frac ba
\sum _{n=0}^\infty \frac{(q^{n+1};q)_\infty (bq^{n+1}/a;q)_\infty}
{(bq^{n+1};q)_\infty} q^n
P_m( bq^{n+1} ;a,b;q)P_{m'}( bq^{n+1} ;a,b;q)
$$  $$
=\frac{(q,b/a,aq/b;q)_\infty}{(aq,bq;q)_\infty} \frac{(q;q)_m}
{(aq,bq;q)_m}(-ab)^{m}q^{m(m+3)/2}\delta_{mm'} . \eqno (3.5)
$$

Notice that for $m=m'=0$ the orthogonality relation (3.5) reduces
to
$$
\sum_{n=0}^\infty \frac{(q^{n+1},aq^{n+1}/b;q)_\infty}
{(aq^{n+1};q)_\infty} q^n-\frac ba
\sum_{n=0}^\infty \frac{(q^{n+1},bq^{n+1}/a;q)_\infty}
{(bq^{n+1};q)_\infty} q^n =
\frac{(q,aq/b,b/a;q)_\infty}{(aq,bq;q)_\infty} .
\eqno (3.6)
$$
In terms of the ${}_2\phi_1$ basic
hypergeometric series this identity can be written as
$$
\frac{(aq/b,q;q)_\infty}{(aq;q)_\infty}\; {}_2\phi_1 (aq,0;\; aq/b;\;
q,q)
$$   $$
- \frac ba
\frac{(bq/a,q;q)_\infty}{(bq;q)_\infty}\; {}_2\phi_1 (bq,0;\; bq/a;\;
q,q) =\frac{(aq/b,b/a,,q;q)_\infty}{(aq,bq;q)_\infty}
\eqno (3.7)
$$
and it represents a particular case of Sears' three-term transformation
formula for ${}_2\phi_1$ series (see [8], formula (3.3.5)).
So the orthogonality relation (3.5) is based on the relation
(3.7). A detailed derivation of this transformation formula and
the orthogonality relation (3.5) can be found in [13].

The explicit form of the orthogonality relation (3.5) is
important for the case under discussion because it directly
leads to the following statement.
\medskip

\noindent {\bf Theorem 1.} {\it The spectrum of the operator $A$
coincides with the set of points $aq^{n+1}$, $bq^{n+1}$,
$n=0,1,2,\cdots$. The spectrum is simple and it has only one
accumulation point at 0.}
\bigskip

\noindent{\bf 4. Representations $T^+_l$, related to big
$q$-Laguerre polynomials}
\bigskip

The operator $A$ has eigenfunctions $\xi_\lambda(x)$, corresponding
to the eigenvalues $aq^{n+1}$, $bq^{n+1}$, $n=0,1,2,\cdots$. Since
these eigenvalues are distinct, the set of  functions
$$
\Xi_n(x) \equiv \xi_{aq^{n+1}}(x),\ \
\Xi'_n(x) \equiv \xi_{bq^{n+1}}(x) \ \ \ \
n=0,1,2,\cdots ,
$$
constitutes an orthogonal basis in the representation space ${\cal
H}_l$. For convenience, we often denote this basis as
$$
\tilde\Xi_n(x)  \ \ \ \ n=0,\pm 1,\pm 2,\cdots ,
$$
where $\tilde\Xi_{n}(x):=\Xi'_{n}(x)$, $n=0,1,2,\cdots$, and
$\tilde\Xi_{-m}(x):=\Xi'_{m-1}(x)$,
$m=1,2,3,\cdots$. This basis is not orthonormal. Let us find the
orthonormal basis
$$
\hat\Xi_n(x)=c_n  \tilde\Xi_n(x),\ \ \ \ \  n=0,\pm 1,\pm 2,\cdots
, \eqno (4.1)
$$
where $c_n$ are normalization constants.

Observe that due to (3.4) the bases $f^l_n$, $n=0,1,2,\cdots$, and
$\hat\Xi_n(x)$, $n=0,\pm 1,\pm 2,\cdots$, are connected by the
formulas
$$
\hat\Xi_n(x)=c_n\sum_{m=0}^\infty a_m(aq^{n+1})f^l_m(x),\ \ \
\hat\Xi'_n(x)=c'_n\sum_{m=0}^\infty a_m(bq^{n+1})f^l_m(x),\ \ \
n=0,1,2,\cdots ,
$$
where $a_m(\lambda)$ are defined in (3.3) and $\hat\Xi'_{n-1}(x)
=\hat\Xi_{-n}(x)$ for $c'_{n-1}=c_{-n}$, $n=1,2,\cdots$. In order
to find the coefficients $c_n$ we take into account the following.
The basis $\hat\Xi_n(x)$, $n=0,\pm 1,\pm 2,\cdots$, is orthonormal
if the matrix $(u_{mn})$ with entries $u_{mn}=c_na_m(\lambda_n)$,
where $\lambda_n=aq^{n+1}$ for $n=0,1,2,\cdots$ and
$\lambda_n=bq^{-n}$ for $n=-1,-2,\cdots$, is unitary, that is,
$\sum _n u_{mn}u_{m'n}=\delta _{mm'}$. Taking into account the
explicit form (3.3) of the coefficients $a_m(\lambda _n)$ and the
orthogonality relation (3.5) for the big $q$-Laguerre polynomials,
we find that
$$
c_n=\left( \frac{(q^{n+1},aq^{n+1}/b,aq,bq;q)_\infty \, q^n}
{(aq^{n+1},q,b/a,aq/b;q)_\infty}\right )^{\frac 12} =\left(
\frac{(aq;q)_n(bq;q)_\infty \, q^n} { (aq/b,q;q)_n
(b/a;q)_\infty}\right )^{\frac 12} \eqno (4.2)
$$
and
$$
c'_n=\left( \frac{(-b/a)q^{n}(bq^{n+1}/a,q^{n+1},aq,bq;q)_\infty}
{(bq^{n+1},q,b/a,aq/b;q)_\infty}\right )^{\frac 12}
=\left( \frac{(-b/a)q^{n} (bq;q)_n(aq;q)_\infty }
{(q;q)_n (aq/b;q)_\infty (b/a;q)_{n+1}}\right )^{\frac 12}.
\eqno (4.3)
$$
Namely, at these values of $c_n$ and $c'_n$ the formula $\sum _n
u_{mn}u_{m'n}=\delta _{mm'}$ is equivalent to the orthogonality
relation for the big $q$-Laguerre polynomials.

We know that the operator $A$ acts upon the basis $\Xi_n(x)$, $\Xi'_n(x)$,
$n=0,1,2,\cdots$, in the following way:
$$
A\, \Xi_n(x)=aq^{n+1}\Xi_n(x),\ \ \ A\,
\Xi'_n(x)=bq^{n+1}\Xi'_n(x),
$$
Let us find how the operator $q^{-J_0}$ acts upon this basis.
From the $q$-difference equation (3.11.5) in [11], it follows that
$$
q^{-n}\lambda^2P_n(\lambda)=B(\lambda) P_n(q\lambda)-[B(\lambda)
+D(\lambda)-\lambda^2]P_n(\lambda)+D(\lambda)P_n(q^{-1}\lambda) ,
$$
where $P_n(\lambda):=P_n(\lambda ;a,b;q)$,
$B(\lambda)=abq(1-\lambda)$, and $D(\lambda)=(\lambda-aq)(\lambda
-bq)$. We multiply both sides of this relation by $b_nf^l_n(x)$,
where $b_n$ is the constant factor in front of the $P_n(\lambda ;
a,b;q)$ on the right-hand side of (3.3), and sum up over $n$.
Taking into account formula (3.4) and the relation
$q^{-J_0}f^l_n=q^{-l-n}f^l_n$, we obtain
$$
q^{-J_0+l}\lambda^2\xi_\lambda (x)=B(\lambda)
\xi_{q\lambda}(x)-[B(\lambda) +D(\lambda)-\lambda^2]\xi_\lambda
(x)+D(\lambda)\xi_{q^{-1}\lambda}(x).
$$
Since $B(\lambda)+D(\lambda)-\lambda^2 = abq(1+q)-\lambda
q(ab+a+b)$, we have
$$
q^{-J_0}\Xi_n=a^{-3/2}bq^{-2n-3/2}(1-aq^{n+1})\Xi_{n+1}  -
a^{-3/2}q^{-2n-3/2} [b(1+q)-q^{n+1}(ab+a+b)]\Xi_n
$$  $$
+ \,a^{-3/2}bq^{-2n-1/2}(1-q^n)(1-aq^n/b)\Xi_{n-1} \eqno (4.4)
$$
for $\lambda =aq^{n+1}$, that is, for $\xi_{aq^{n+1}} (x)=
\Xi_n(x)$, and
$$
q^{-J_0}\Xi'_n=a^{1/2}b^{-1}q^{-2n-3/2}(1-bq^{n+1})\Xi'_{n+1}-
a^{1/2}b^{-1}q^{-2n-3/2}
$$  $$
\times [1+q-a^{-1}q^{n+1}(ab+a+b)] \Xi'_n
+a^{1/2}b^{-1}q^{-2n-1/2}(1-q^n)(1-bq^{n}/a)\Xi'_{n-1}  \eqno (4.5)
$$
for $\lambda =bq^{n+1}$, that is, for $\xi_{bq^{n+1}} (x)=
\Xi'_n(x)$. Passing in (4.4) and (4.5) to the orthonormal basis
$\hat\Xi_n$, $n=0,\pm 1,\pm 2,\cdots$, we obtain
$$
q^{-J_0}\hat\Xi_n=a^{-3/2}bq^{-2n-2}\sqrt{(1-aq^{n+1})
(1-q^{n+1})(1-aq^{n+1}/b)}\, \hat\Xi_{n+1}
$$  $$
-a^{-3/2}q^{-2n-3/2} [b(1+q)-q^{n+1}(ab+a+b)]\, \hat\Xi_n
$$  $$
+a^{-3/2}bq^{-2n}\sqrt{(1-aq^{n})(1-q^n)(1-aq^n/b)} \,
\hat\Xi_{n-1}
$$
for $\hat\Xi_n$, $n=0,1,2,\cdots$, and
$$
q^{-J_0}\hat\Xi'_n=a^{1/2}b^{-1}q^{-2n-2}\sqrt{(1-bq^{n+1})
(1-q^{n+1})(1-bq^{n+1}/a)} \, \hat\Xi'_{n+1}
$$  $$
-a^{1/2}b^{-1}q^{-2n-3/2} [1+q-a^{-1}q^{n+1}(ab+a+b)]\, \Xi'_n
$$  $$
+a^{1/2}b^{-1}q^{-2n}\sqrt{(1-bq^{n})(1-q^n)(1-bq^{n}/a)} \,
\hat\Xi'_{n-1}.
$$
The operators $A$ and $q^{-J_0}$ completely determine the
representation $T^+_l$ in the basis $\hat\Xi_n$, $n=0,\pm 1,\pm 2,
\cdots$ (see, for example, [14]). However, the expressions for
the operators $J_+$ and $J_-$ are rather complicated.
\bigskip

\noindent{\bf 5. Dual polynomials and functions}
\bigskip

The matrix $(u_{mn})$, $m=0,1,2,\cdots$, $n=0,\pm 1,\pm 2,\cdots$,
with entries $u_{mn}=c_na_m(\lambda_n)$,
described in the previous section, is unitary and it connects two
orthonormal bases in the Hilbert space ${\cal H}_l$.
The unitarity of this matrix means
that the following relations hold:
$$
\sum_{n\in {\Bbb Z}} u_{mn}u_{m'n}=\delta_{mm'},\ \ \
\sum_{m=0}^\infty u_{mn}u_{mn'}=\delta_{nn'} . \eqno (5.1)
$$
It is easy to see that the first relation is equivalent to
the orthogonality relation for the big $q$-Laguerre polynomials (see
the previous section).
The second relation is the orthogonality relation for the
functions, which are dual to the big $q$-Laguerre
polynomials, and are defined as
$$
f_n(q^{-m}; a,b|q):= P_m(aq^{n+1};a,b;q), \ \ \ n=0,1,2,\cdots ,
\eqno (5.2)
$$  $$
g_n(q^{-m}; a,b|q):= P_m(bq^{n+1};a,b;q), \ \ \ n=0,1,2,\cdots .
\eqno (5.3)
$$
Taking into account the expressions for the entries $u_{mn}$, the
relation $\sum_{m=0}^\infty u_{mn}u_{mn'}=\delta_{nn'}$ can be
written as
$$
\sum _{m=0}^\infty a_m(\lambda_n)a_m(\lambda_{n'})=c_n^{-2}
\delta _{nn'},
$$
where $c_n$ must be replaced by $c'_n$ if $\lambda_n =bq^{n+1}$.
Substituting the explicit expressions for the coefficients
$a_m(\lambda_n)$, we derive the following orthogonality relations
for the functions (5.2) and (5.3):
$$
\sum _{m=0}^\infty \frac{(aq,bq;q)_m}
{(q;q)_m(-abq^{2})^{m}} q^{-m(m-1)/2} f_n(q^{-m};
a,b|q)f_{n'}(q^{-m}; a,b|q) =c_n^{-2}\delta_{nn'}, \eqno (5.4)
$$  $$
\sum _{m=0}^\infty \frac{(aq,bq;q)_m}
{(q;q)_m(-abq^{2})^{m}}q^{-m(m-1)/2}\, g_n(q^{-m};
a,b|q)\, g_{n'}(q^{-m}; a,b|q) ={c'}_n^{-2}\delta_{nn'}, \eqno (5.5)
$$  $$
\sum _{m=0}^\infty \frac{(aq,bq;q)_m}
{(q;q)_m(-abq^{2})^{m}}q^{-m(m-1)/2}\, f_n(q^{-m};
a,b|q)\, g_{n'}(q^{-m}; a,b|q)=0,
 \eqno (5.6)
$$
where $c_n$ and $c'_n$ are given by the formulas (4.2) and
(4.3).

Comparing the expression (3.13.1) in [11] for the $q$-Meixner
polynomials
$$
M_n(q^{-x};a,b;q):={}_2\phi_1 (q^{-n},q^{-x};\; aq;\ q,-q^{n+1}/b)
$$
with the explicit form (3.2) of
the big $q$-Laguerre polynomials $P_m (x;a,b;q)$, we see that
$$
f_n(q^{-m};a,b|q)=(q^{-m}/b;q)^{-1}_m M_n(q^{-m};a,-b/a;q).
$$
Since $(q^{-m}/b;q)_m=(bq;q)_m(-b)^{-m}q^{-m(m+1)/2}$,
the orthogonality relation (5.4) leads to the
orthogonality relation for the $q$-Meixner polynomials
$M_n(q^{-m})\equiv M_n(q^{-m}; a,-b/a;q)$:
$$
\sum _{m=0}^\infty \frac{(aq;q)_m (-b/a)^m q^{m(m-1)/2}
}{(bq,q;q)_m} M_n(q^{-m})M_{n'}(q^{-m})
$$
$$
=\frac{(b/a;q)_\infty}{(bq;q)_\infty}
\frac{(aq/b,q;q)_n}{(aq;q)_n} q^{-n} \delta_{nn'} ,
 \eqno (5.7)
$$
where, as before, $0<a<q^{-1}$ and $b<0$. This orthogonality
relation coincides with formula (3.13.2) in [11].

The functions (5.3) are also expressed in terms of $q$-Meixner
polynomials. Indeed, we have
$$
g_n(q^{-m};a,b|q)={}_3\phi_2 (q^{-m},0,bq^{n+1};\; aq,bq;\ q,q)
$$  $$
=(q^{-m}/a;q)_m^{-1}\, {}_2\phi_1 (q^{-n},q^{-m};\; bq;\
q,bq^{n+1}/a)
$$  $$
=(q^{-m}/a;q)_m^{-1}\, M_n(q^{-m}; b,-a/b;\; q) ,
$$
where $b<0$, that is, one of the parameters in these $q$-Meixner
polynomials is negative.

Substituting this expression for $g_n(q^{-m};a,b|q)$ into (5.5),
we obtain the orthogonality relation for $q$-Meixner polynomials
$M_n(q^{-m})\equiv M_n(q^{-m}; b,-a/b;\; q)$ with negative $b$:
$$
\sum _{m=0}^\infty \frac{(bq;q)_m (-a/b)^m}{(aq,q;q)_m}
\, q^{m(m-1)/2}\,M_n(q^{-m})M_{n'}(q^{-m})
$$
$$
=\frac{(a/b;q)_\infty}{(aq;q)_\infty}
\frac{(bq/a,q;q)_n}{(bq;q)_n} q^{-n}\delta_{nn'}. \eqno (5.8)
$$
Observe that this orthogonality relation is of the same form as
for $b>0$ (see, for example, formula (3.13.2) in [11]). As far as
we know, this type of orthogonality relation for negative values
of the parameter $b$ has not been discussed in the literature.

The relation (5.6) can be written as the equality
$$
\sum_{m=0}^\infty \frac{(-1)^m q^{m(m-1)/2}}{(q;q)_m}
 M_n(q^{-m}; a,-b/a;\; q)M_{n'}(q^{-m}; b,-a/b;\; q)=0,
 $$
which holds for $n,n'=0,1,2,\cdots$. The reader is invited to verify
directly the validity of this identity for arbitrary nonnegative
integers $n$ and $n'$ by using Jackson's $q$-exponential function
$$
E_q(z):= \sum_{n=0}^\infty \frac{q^{n(n-1)/2}}{(q;q)_n}z^n=
(-z;q)_\infty
$$
and the fact that $E_q(z)$ has zeroes at the points $z_j = -q^{-j}$,
$j=0,1,2,\cdots$, namely, $E_q(-q^{-j})=0$.

Notice that the appearance of the $q$-Meixner polynomials
here as a dual family with respect to the big $q$-Laguerre
polynomials is quite natural because the transformation
$q\to q^{-1}$ interrelates these two sets of polynomials,
that is,
$$
M_n(x;b,c;q^{-1})=(q^{-n}/b;q)_n\, P_n(qx/b;1/b,-c;q).
$$

Let us introduce the Hilbert space ${\frak l}_b^2$ of functions
$F(q^{-m})$ on the set $m\in \{0,1,2,\cdots \}$ with a scalar
product given by the formula
$$
\langle F_1,F_2\rangle_b = \sum _{m=0}^\infty \frac{(aq,bq;q)_m}
{(q;q)_m(-abq^2)^{m}}\,
q^{-m(m-1)/2} \,F_1(q^{-m})\overline{F_2(q^{-m})}. \eqno (5.9)
$$
Now we can formulate the following statement.
\medskip

{\bf Theorem 2.} {\it The functions (5.2) and (5.3) constitute an
orthogonal basis in the Hilbert space ${\frak l}_b^2$.}
\medskip

{\it Proof.} To show that the system of functions (5.2) and (5.3)
constitutes a complete basis in the space ${\frak l}_b^2$ we take
in ${\frak l}_b^2$ the set of functions $F_k$, $k=0,1,2, \cdots$,
such that $F_k(q^{-m})=\delta _{km}$. It is clear that these
functions constitute a basis in the space ${\frak l}_b^2$. Let us
show that each of these functions $F_k$ belongs to the closure
$\bar V$ of the linear span $V$ of the functions (5.2) and (5.3).
This will prove the theorem. We consider the functions
$$
{\hat F}_k(q^{-m})=\sum _{n=-\infty} ^\infty u_{kn}u_{mn},\ \ \ \
k=0,1,2,\cdots ,
$$
where $u_{jn}$ are the same as in (5.1). Then $\hat F_k(q^{-m})\in
\bar V$ and, due to the first equality in (5.1), $\hat F_k$,
$k=0,1,2,\cdots$, coincide with the corresponding functions $F_k$,
introduced above. The theorem is proved.
\medskip

The measure in (5.9) does not coincide with the orthogonality
measure for $q$-Meixner polynomials. Multiplying the measure in
(5.9) by $[(bq;q)_m(-b)^{-m}q^{-m(m+1)/2}]^{-2}$, we obtain the
measure in (5.7). Let ${\frak l}^2_{(1)}$ be the Hilbert space of
functions $F(q^{-m})$ on the set $m\in \{0,1,2,\cdots \}$ with the
scalar product
$$
\langle F_1,F_2\rangle_{(1)} = \sum _{m=0}^\infty \frac{(aq;q)_m
(-b/a)^m q^{m(m-1)/2}} {(bq,q;q)_m}\,F_1(q^{-m})\,\overline{F_2(q^{-m})} ,
$$
where the weight function coincides with the measure in (5.7).

Taking into account the modification of the measure and the
statement of Theorem 2, we conclude that the $q$-Meixner
polynomials $M_n(q^{-m};\, a,-b/a;\; q)$ and the functions
$$
(bq;q)_m(-b)^{-m}q^{-m(m+1)/2} g_n(q^{-m}; a,b|q)
$$
constitute an orthogonal basis in the space ${\frak l}^2_{(1)}$.
\medskip

{\bf Proposition 1.} {\it The $q$-Meixner polynomials
$M_n(q^{-m};a,c;q)$, $n=0,1,2,\cdots$, with the parameters
$a=q^{2l-1}$ and $c=-b/a$ do not constitute a complete basis in
the Hilbert space ${\frak l}^2_{(1)}$, that is, the $q$-Meixner
polynomials are associated with the indeterminate moment problem
and the measure in (5.7) is not an extremal measure for these
polynomials.}
\medskip

{\it Proof.} In order to prove this proposition we note that if
the $q$-Meixner polynomials would be associated with the determinate
moment problem, then they would constitute a basis in the space of
squared integrable functions with respect to the measure from
(5.7). However, this is not the case. By the definition of an
extremal measure, if the measure in (5.7) would be extremal (see,
for example, [12], Chapter VII), then again the set of the
$q$-Meixner polynomials would be a basis in that space. Therefore,
the measure is not extremal. Proposition is proved.
\medskip

Let now ${\frak l}^2_{(2)}$ be the Hilbert space of functions
$F(q^{-m})$ on the set $m\in \{0,1,2,\cdots \}$, with the scalar
product
$$
\langle F_1,F_2\rangle_{(2)} = \sum _{m=0}^\infty \frac{(bq;q)_m
(-a/b)^m q^{m(m-1)/2}} {(aq,q;q)_m}\,F_1(q^{-m})\,\overline{F_2(q^{-m})}.
$$
The measure above coincides with the orthogonality measure in (5.8)
for $q$-Meixner polynomials $M_n(q^{-m}; b,-a/b;\; q)$, $b<0$. The
following proposition is proved in the same way as Proposition 2.
\medskip

{\bf Proposition 2.} {\it The $q$-Meixner polynomials
$M_n(q^{-m};b,-a/b;q)$, $n=0,1,2,\cdots$, with $b<0$ do not
constitute a complete basis in the Hilbert space ${\frak
l}^2_{(2)}$, that is, these $q$-Meixner polynomials are
associated with the indeterminate moment problem and the
measure in (5.8) is not an extremal measure for them.}
\medskip

According to Propositions 1 and 2, the measures in formulas (5.7)
and (5.8), with respect to which the $q$-Meixner polynomials are
orthogonal, are not extremal. As far as we know, explicit forms of
extremal measures for the $q$-Meixner polynomials are not known.
It is worth to mention that extremal measures have been constructed
only for the $q$-Hermite polynomials when $q>1$ (see [15]).

Note that the complementary set of orthogonal functions to the
$q$-Laguerre polynomials $L_n^{(\alpha )}(x;\,q)$ in the Hilbert
space of square integrable functions with respect to the orthogonality
measures (which are not extremal) for these polynomials has been
constructed in [16].
\bigskip

\noindent{\bf 6. Generating function for big $q$-Laguerre polynomials}
\bigskip

The aim of this section is to derive a generating function
for the big $q$-Laguerre polynomials
$$
G(x,t; a,b;\, q):= \sum _{n=0}^\infty
\frac{(aq,bq;q)_nq^{-n(n-1)/2}}{(q;q)_n}P_n(x;a,b;q)t^n ,
\eqno (6.1)
$$
which will be used in the next section. Observe that this
formula (6.1) is a bit more general than each of the three
instances of generating functions for big $q$-Laguerre
polynomials, given in section 3.11 of [11].

Employing the explicit expression
$$
P_n(x;a,b;q)=(b^{-1}q^{-n};q)_n^{-1}{}_2\phi_1(q^{-n},aqx^{-1};\;
aq;\ q, x/q)
$$
for the big $q$-Laguerre polynomials, one obtains
$$
G(x,t; a,b;\, q)=\sum_{n=0}^\infty \frac{(aq;q)_n}{(q;q)_n}
(-bqt)^n \sum_ {k=0}^n \frac{(q^{-n},aqx^{-1};q)_k}{(aq,q;q)_k}
\left( \frac xb \right) ^k
$$
$$
=\sum_{n=0}^\infty (aq;q)_n
(-bqt)^n \sum_ {k=0}^n \frac{(-x/b)^k(aqx^{-1};q)_k}
{(aq,q;q)_k(q;q)_{n-k}} q^{-nk+k(k-1)/2}
$$  $$
=\sum_{k=0}^\infty \frac{(aqx^{-1};q)_k(-x/b)^k}
{(aq,q;q)_k} q^{k(k-1)/2}
 \sum_ {m=0}^\infty \frac{(aq;q)_{m+k}}
{(q;q)_m} (-bqt)^{m+k} q^{-(k+m)k}
$$   $$
=\sum_{k=0}^\infty \frac{(aqx^{-1};q)_k}
{(q;q)_k} (xt)^k q^{-k(k-1)/2}
 \sum_ {m=0}^\infty \frac{(aq^{k+1};q)_m}
{(q;q)_m} (-bq^{1-k}t)^m .
$$
By the $q$-binomial theorem, the last sum equals to
$(-abq^2;q)_\infty/(-bq^{1-k}t;q)_\infty$.
Since
$$
(-bq^{1-k}t;q)_\infty =q^{-k(k-1)/2}(-q/bqt;q)_k(-bqt;q)_\infty ,
$$
then
$$
\frac{(-abq^2;q)_\infty}{(-bq^{1-k}t;q)_\infty}
=\frac{(-abq^2;q)_\infty}{(-bqt;q)_\infty}
\frac{q^{k(k-1)/2}}{(bt)^k(-1/bt;q)_k}.
$$
Thus,
$$
G(x,t; a,b;\, q)=\frac{(-abq^2;q)_\infty}{(-bqt;q)_\infty}
\sum_{k=0}^\infty \frac{(aqx^{-1};q)_k} {(-1/bt,q;q)_k}\left(
\frac xb \right) ^k
$$   $$
 =\frac{(-abq^2;q)_\infty}{(-bqt;q)_\infty}
{}_2\phi_1 (aqx^{-1},\, 0;\, -1/bt;\; q,x/b). \eqno (6.2)
$$
This gives a desired generating function for the big
$q$-Laguerre polynomials.
\bigskip
\bigskip

\noindent{\bf 7. Biorthogonal systems of functions}
\bigskip

From the very beginning we could consider an operator
$$
A_1 :=\alpha q^{J_0/4}\,\left[(1-bq^{J_0-l})J_+ +q^{J_0-l} J_- \right]\,
 q^{J_0/4} - \beta_1 q^{2J_0} + \beta_2 q^{J_0-l}
$$
(cf (3.1)), where $\alpha, \beta_1$, and $\beta_2$ are the same as
in (3.1). This operator is well defined, but it is not selfadjoint.
Repeating the reasoning of section 3, we find that eigenfunctions
of $A_1$ are of the form
$$
\psi_\lambda(x)=
\sum _{n=0}^\infty
(-ab)^{-n/2}q^{-n}\left( \frac{(aq;q)_n
}{(q;q)_n} \right) ^{1/2}
P_n( \lambda ;a,b;q) f^l_n (x)
$$  $$
=\sum _{n=0}^\infty a^{-3n/4}(-b)^{-n/2}q^{-n} \frac{(aq;q)_n}
{(q;q)_n} P_n( \lambda ;a,b;q)x^n , \eqno (7.1)
$$
where, as before, $a=q^{2l-1}$.
The last sum can be summed with the aid of formula (3.11.12)
in [11]. We thus have
$$
\psi_\lambda(x)=((-a/b^2)^{1/4}x;q)_\infty \cdot {}_2\phi_1
(bq\lambda^{-1},\, 0;\, bq;\ q, a^{-3n/4}(-b)^{-1/2}q^{-1}x\lambda
).
$$

Now we consider another operator
$$
A_2 :=\alpha q^{J_0/4}(J_+q^{J_0-l} + J_-
(1-bq^{J_0-l})) q^{J_0/4}
-\beta_1 q^{2J_0} +\beta_2 q^{J_0-l},
$$
where $\alpha, \beta_1$, and $\beta_2$ are
the same as above. This operator is adjoint to the operator $A_1:
A_2^*=A_1$. Repeating the reasoning of section 3, we find that
 eigenfunctions of $A_2$ have the form
$$
\varphi_\lambda(x)=
\sum _{n=0}^\infty
(-ab)^{-n/2}q^{-n(n+1)/2}\left( \frac{(aq;q)_n
(bq;q)^2_n}{(q;q)_n} \right) ^{1/2}
P_n( \lambda ;a,b;q) f^l_n (x)
$$  $$
=\sum _{n=0}^\infty
a^{-3n/4}(-b)^{-n/2}q^{-n(n+1)/2} \frac{(aq;q)_n (bq;q)_n}
{(q;q)_n}  P_n( \lambda ;a,b;q)x^n. \eqno (7.2)
$$
According to the formula (6.2), this function
can be written as
$$
\varphi_\lambda(x)=\frac{(-abq^{2};q)_\infty}
{(a^{-3/4}(-b)^{1/2}x;q)_\infty}\, {}_2\phi_1(
aq/\lambda ,0;\, a^{3/4}(-b)^{-1/2}q/x;\;
q,\lambda /b).
$$

Let us denote by $\Psi_m(x)$, $m=0,\pm 1,\pm 2,\cdots$, the functions
$$
\Psi_m(x)=c_m\psi_{aq^{m+1}}(x),\ \ m=0,1,2,\cdots,\ \ \
\Psi_{-m}(x)=c'_{m-1}\psi_{bq^m}(x),\ \ m=1,2,\cdots,
 \eqno (7.3)
$$
and by $\Phi_m(x)$, $m=0,\pm 1,\pm 2,\cdots$, the functions
$$
\Phi_m(x)=c_m\varphi_{aq^{m+1}}(x), \ \ m=0,1,2,\cdots,\ \ \
\Phi_{-m}(x)=c'_{m-1}\varphi_{bq^{m}}(x), \ \ m=1,2,\cdots,
 \eqno (7.4)
$$
where $c_m$ and $c'_m$ are given by formulas (4.2) and (4.3).

Writing down the decompositions (7.1) and (7.2) for the functions
$\Psi_m(x)$ and $\Phi_m(x)$ (in terms of the orthonormal
basis $f^l_n$, $n=0,1,2,\cdots$, of the Hilbert space ${\cal H}_l$)
and taking into account the orthogonality relations (5.4)--(5.6)
we find that
$$
\langle \Psi_m(x),\Phi_n(x)\rangle =\delta_{mn}, \ \ \
m,n=0,\pm 1,\pm 2,\cdots .
$$
This means that we can formulate the following statement.
\medskip

{\bf Theorem 3.} {\it
The set of functions $\Psi_m(x)$, $m=0,\pm 1,\pm 2,\cdots$,
and the set of functions $\Phi_m (x)$, $m=0,\pm 1,\pm 2,\cdots$,
form biorthogonal sets of functions with respect to the
scalar product in the Hilbert space ${\cal H}_l$.}
\bigskip

\noindent{\bf 8. The classical limit as $q\to 1$}
\bigskip

In section 3 we have shown that the operator $A$, defined by
(3.1), is related to the family of big $q$-Laguerre polynomials
(3.2). Namely, the eigenfunctions $\xi_\lambda(x)$ of the operator
$A$ can be expanded in terms of the canonical basis functions of
the representation $T_l^+$ of the algebra $U_q({\rm su}_{1,1})$
and coefficients of this expansion are big $q$-Laguerre
polynomials (up to multiplication by a constant factor, see
(3.4)).

It is well known that in the limit as $q\to 1$ big $q$-Laguerre
polynomials $P_n(x;a,b;q)$ reduce to classical Laguerre
polynomials $L_n^{(\alpha)}(x)$, i.e.,
$$
{\rm lim}_{q\to 1} P_n(x;q^\alpha ,(q-1)^{-1}q^\beta ;q)
=\frac{L_n^{(\alpha)}(1-x)}{L_n^{(\alpha)}(0)}. \eqno (8.1)
$$
So, it is natural to expect that there exists a corresponding
classical limit of the operator $A$ as $q\to 1$ and $b=q^\beta
/(q-1)\to -\infty$. Indeed, this is the case. Bearing in mind that
in the case under discussion $a=q^{2l-1}$ and substituting
$b=q^\beta /(q-1)$ into (3.1), it is not hard to evaluate that
$$
A^{\rm cl}:= {\rm lim}_{q\to 1}A=2(J_1^{\rm cl}-J_0^{\rm cl})+I,
\eqno (8.2)
$$
where $I$ is the identity operator and $J_1^{\rm cl}$ and
$J_0^{\rm cl}$ are the generators of the classical Lie algebra
${\rm su}(1,1)$, explicitly realized in terms of the first-order
differential operators:
$$
J_0^{\rm cl}=x\frac{d}{dx} +l,\ \ \ \ J_1^{\rm cl}=\frac 12
(1+x^2)\frac {x}{dx} +lx. \eqno (8.3)
$$

Observe that in this particular limit the coefficients of the
expansion of the eigenfunctions $\xi_\lambda (x)$ in monomials
$x^n$ (see the second line in (3.4)) tend to the Laguerre
polynomials $L_n^{(2l-1)}(1-\lambda)$. This means that
$$
\xi^{\rm cl}_\lambda (x):= {\rm lim}_{q\to 1}\xi_\lambda (x)
=\sum_{n=0}^\infty L_n^{(2l-1)}(1-\lambda)x^n =
\frac{1}{(1-x)^{2l}} \, \exp \, \left( \frac{(\lambda-1)x}{1-x}
\right) , \eqno (8.4)
$$
where at the last step we employed the generating function
$$
\sum_{n=0}^\infty L_n^{(\alpha)}(x)t^n = \frac{1}{(1-t)^{\alpha+1}}
\,\exp \left(\frac{xt}{t-1} \right)
$$
for the classical Laguerre polynomials.

As a consistency check, it is now easy to evaluate directly with
the aid of formulas (8.2)--(8.4) that $A^{\rm cl} \xi^{\rm
cl}_\lambda (x)=\lambda\xi^{\rm cl}_\lambda (x)$.
\bigskip

\noindent{\bf 9. Concluding remarks}
\bigskip

It is well known that many physical systems admit symmetries with
respect to the quantum algebra $U_q({\rm su}_{1,1})$. Then the
collection of states of such a system forms a representation space
for this quantum algebra (see, for example, [17]). Hamiltonians of
such systems often coincide with representation operators of
$U_q({\rm su}_{1,1})$. Since spectra of Hamiltonians are bounded
from below, these representations belong, as a rule, to the
positive discrete series.

Representation operators of quantum algebras are much more
complicated than in the case of Lie algebras. For example, many
symmetric representation operators are unbounded and their
closures are not selfadjoint operators. Such operators have
selfadjoint extensions only if their deficiency indices are equal
to each other. In the case of representations of the algebra
$U_q({\rm su}_{1,1})$ many representation operators can be
equivalently written in the form of a Jacobi matrix. In this case,
the selfadjointness of representation operators causes quite
definite properties of the corresponding polynomial families. In
other words, the study of representation operators leads to deeper
understanding of orthogonal polynomials.

In the present paper we have studied in detail those operators in
the discrete series representations of $U_q({\rm su}_{1,1})$,
which are associated with big $q$-Laguerre polynomials. In
particular, we have found an explicit form of the spectrum for the
corresponding selfadjoint representation operator and derived
explicitly its eigenfunctions. Using these results, we constructed
a system of orthogonal polynomials dual to the big $q$-Laguerre
polynomials. It has occurred that this dual system consists of two
sets of $q$-Meixner polynomials. It is deduced from this fact that
the family of $q$-Meixner polynomials with fixed parameters does
not constitute a complete basis in the $L^2$ space with respect to
their orthogonality measure. This means that $q$-Meixner
polynomials correspond to a representation operator with a
closure, which is not selfadjoint and has deficiency indices
(1,1). Selfadjoint extensions of this operator correspond to so
called extremal measures of orthogonality for $q$-Meixner
polynomials. (The knowledge of these measures would give us a
possibility to find explicitly spectra of selfadjoint extensions
of the operator). Unfortunately, these measures are not known. We
hope that the further development of the approach of this paper
will enable us to find extremal measures for orthogonal
polynomials, which correspond to an indeterminate moment problem.
 \bigskip

\noindent{\bf Acknowledgments}

\medskip

This research has been supported in part by the SEP-CONACYT
project 41051-F and the DGAPA-UNAM project IN112300 "Optica
Matem\'atica". A. U. Klimyk acknowledges the Consejo Nacional de
Ciencia y Technolog\'{\i}a (M\'exico) for a C\'atedra Patrimonial
Nivel II.

\bigskip

\noindent{\bf References}
\bigskip

1. T. H. Koornwinder, {\it Askey--Wilson polynomials as zonal
spherical functions on the $SU(2)$ quantum group}, SIAM J. Math.
Anal. {\bf 24}, 795--813 (1993).

2. M. Noumi and K. Mimachi,  {\it Askey--Wilson polynomials and
the quantum group $SU_q(2)$}, Proc. Japan Acad. Ser. A Math. {\bf
66}, 146--149 (1990).

3. H. T. Koelink, {\it The addition formula for continuous
$q$-Legendre polynomials and associated spherical elements on the
$SU(2)$ quantum group related to Askey--Wilson polynomials}, SIAM
J. Math. Anal. {\bf 25}, 197--217 (1994).

4. H. T. Koelink and J. Van der Jeugt, {\it Convolutions for
orthogonal polynomials from Lie and quantum algebra
representations}, SIAM J. Math. Anal. {\bf 29} (1998), 794--822.

5. J. Van der Jeugt and R. Jagannathan, {\it Realizations of
$su(1,1)$ and $U_q(su(1,1))$ and generating functions for
orthogonal polynomials}, J. Math. Phys. {\bf 39} (1998),
5062--5078.

6. H. T. Koelink and J. Van der Jeugt, {\it Bilinear generating
functions for orthogonal polynomials}, Constructive Approximation
{\bf 14} (1999), 481--497.

7. N. M. Atakishiyev and A. U. Klimyk, {\it Diagonalization of
representation operators for the quantum algebra $U_q({\rm
su}_{1,1})$}, Methods of Functional Analysis and Topology {\bf 8},
No. 3 (2002), 1--12.

8. G. Gasper and M. Rahman, {\it Basic Hypergeometric Functions},
Cambridge University Press, Cambridge, 1990.

9. G. E. Andrews, R. Askey, and R. Roy, {\it Special Functions},
Cambridge  University Press, Cambridge, 1999.

10. I. M. Burban and A. U. Klimyk, {\it Representations of the
quantum algebra $U_q({\rm su}_{1,1})$}, J. Phys. A: Math. Gen.
{\bf 26} (1993), 2139--2151.

11. R. Koekoek and R. F. Swarttouw, {\it The Askey-Scheme of
Hypergeometric Orthogonal Polynomials and Its $q$-Analogue}, Delft
University of Technology Report 98--17; available from {\tt
ftp.tudelft.nl}.

12. Ju. M. Berezanskii, {\it Expansions in Eigenfunctions of
Selfadjoint Operators}, Providence, R. I., Amer. Math. Soc., 1968.

13. M. E. R. Ismail and C. A. Libis, {\it Contiguous relations,
basic hypergeometric functions, and orthogonal polynomials.I}, J.
Math. Anal. Appl. {\bf 141} (1989), 349--372.

14. N. M. Atakishiyev and A. U. Klimyk, {\it Diagonalization of
operators and one-parameter families of nonstandard bases for
representations of ${\rm su}_q(2)$}, J. Phys. A: Math. Gen. {\bf
35} (2002), 5267--5278.

15. M. E. R. Ismail and D. R. Masson, {\it $q$-Hermite
polynomials, biorthogonal rational functions and $q$-beta
integral}, Trans. Amer. Math. Soc. {\bf 346} (1994), 63--116.

16. N. Cicconi, E. Koelink, and T. H. Koornwinder, {\it
$q$-Laguerre polynomials and big $q$-Bessel functions and their
orthogonality relations}, Methods of Applied Analysis {\bf 6}
(1999), 109--127.

17. A. Ballesteros and S. M. Chumakov, {\it On the spectrum of a
Hamiltonian defined on ${\rm su}_q(2)$ and quantum optical
models}, J. Phys. A: Math. Gen. {\bf 32} (1999), 6261--6269.

\end{document}